\documentclass[10pt]{bmc_article}

\usepackage{cite} 
\usepackage{url}  
\usepackage{ifthen}  
\usepackage{multicol}   
\usepackage[utf8]{inputenc} 
\usepackage{hyperref}
\usepackage{amsfonts}
\usepackage{amsmath}
\usepackage{graphicx}
\usepackage{multirow} 
\usepackage{rotating}

\usepackage{ulem}

\newboolean{publ}

\newenvironment{bmcformat}{\begin{raggedright}\baselineskip20pt\sloppy\setboolean{publ}{false}}{\end{raggedright}\baselineskip20pt\sloppy}

\newcommand{\hl}{\hi}
\newcommand{\comment}[1]{} 

\begin{document}
\begin{bmcformat}


\newenvironment{example}[1][Example]{\begin{trivlist}
\item[\hskip \labelsep {\bfseries #1}]}{\end{trivlist}}

\title{ADAM: Analysis of Discrete Models of Biological
Systems Using Computer Algebra}
 

\author{Franziska Hinkelmann$^{1,2}$%
       \email{Franziska Hinkelmann - fhinkel@vt.edu}%
      \and
         Madison Brandon$^{3,\dagger}$%
         \email{Madison Brandon - mbrando1@utk.edu}
		\and
         Bonny Guang$^{4,\dagger}$%
         \email{Bonny Guang - bonny.guang@gmail.com}
		\and
	     Rustin McNeill$^{5,\dagger}$%
	      \email{Rustin McNeill - rcmcneil@uncg.edu}
		\and
         Grigoriy Blekherman$^1$%
          \email{Grigoriy Blekherman - grrigg@vbi.vt.edu}	
		\and
         Alan Veliz-Cuba$^{1,2}$%
         \email{Alan Veliz-Cuba - alanavc@vt.edu}
       and 
         Reinhard Laubenbacher\correspondingauthor$^{1,2}$%
         \email{Reinhard Laubenbacher\correspondingauthor - reinhard@vbi.vt.edu}%
      }


\address{%
    \iid(1)Virginia Bioinformatics Institute, Blacksburg, VA 24061-0123, USA\\
    \iid(2)Virginia Polytechnic Institute and State University, Blacksburg, VA 24061-0123, USA\\
    \iid(3)University of Tennessee - Knoxville, Knoxville, TN 37996-2513, USA\\
    \iid(4)Harvey Mudd College, Claremont, CA 91711-5901, USA\\
    \iid(5)University of North Carolina - Greensboro, Greensboro, NC 27402-6170, USA\\
    \iid(\dagger)These authors contributed equally
}%

\maketitle


\begin{abstract}
        \paragraph*{Background:} Many biological systems are modeled qualitatively with discrete models, such as probabilistic Boolean networks, logical models, Petri nets, and agent-based models, \hl{with the goal to gain a better understanding of the system. The computational complexity to analyze the complete dynamics of these models grows exponentially in the number of variables, which impedes working with complex models. Although there exist sophisticated algorithms to determine the dynamics of discrete models, their implementations usually require labor-intensive formatting of the model formulation, and they are oftentimes not accessible to users without programming skills. Efficient analysis methods are needed that are accessible to modelers and easy to use.}
      
        \paragraph*{Method:} \hl{By converting discrete models into algebraic models, tools from computational algebra can be used to analyze their dynamics. Specifically, we propose a method to identify attractors of a discrete model that is equivalent to solving a system of polynomial equations, a long-studied problem in computer algebra.}

		\paragraph*{Results:} \hl{A method for efficiently identifying attractors, and the web-based tool Analysis of Dynamic Algebraic Models (ADAM), which provides this and other analysis methods for discrete models. ADAM converts several discrete model types automatically into polynomial dynamical systems and analyzes their dynamics using tools from computer algebra. Based on extensive experimentation with both discrete models arising in systems biology and randomly generated networks, we found that the algebraic algorithms presented in this manuscript are fast for systems with the structure maintained by most biological systems, namely sparseness, i.e., while the number of nodes in a biological network may be quite large, each node is affected only by a small number of other nodes, and robustness, i.e., small number of attractors. For a large set of published complex discrete models, ADAM identified the attractors in less than one second.}
		
		\paragraph*{Conclusion:} \hl{Discrete modeling techniques are a useful tool for analyzing complex biological systems and there is a need in the biological community for accessible efficient analysis tools. ADAM provides analysis methods based on mathematical algorithms as a web-based tool for several different input formats, and it makes analysis of complex models accessible to a larger community, as it is platform independent as a web-service and does not require understanding of the underlying mathematics.}
\end{abstract}

\ifthenelse{\boolean{publ}}{\begin{multicols}{2}}{}


\section*{Background}
Mathematical modeling is a crucial tool in understanding the dynamic behavior
of complex biological systems. Discrete models are now widely used for this
purpose. Model types include (probabilistic) Boolean networks, logical
networks, Petri nets, cellular automata, and agent-based (individual-based)
models, to name the most commonly found ones \cite{Steggles, heiner-petri,
shmulevich, Chaouiya, cell-automata, springer_book}. 

There are several existing software packages for simulation and analysis of discrete models. \hl{These include} GINsim, BoolNet R, Snoopy, \hl{signaling Petri net-based simulator in the PathwayOracle toolkit}, DDLab, \hl{GenYsis-P Toolbox }, and BN/PBN Matlab Toolbox \cite{GINsim,boolnet,Snoopy,DDLab,shmulevich-matlab,Garg2008,Dubrova2010,Ruths2008}. Each package has been developed to suit the needs of a particular community, and each package is designed to analyze a different model type. They are discussed in detail in {\it Results and Discussion}. 

Simulation or exhaustive enumeration of the state space is common practice to analyze discrete models, but it is limited by computational complexity, as the state spaces grows exponentially in the number of variables. \hl{When relying on simulation on a standard desktop computer, most software tools are limited to 40 variables or less, i.e., $10^{13}$ states for Boolean systems. Complex models can be analyzed with SAT or BDD-based (Binary Decision Diagram) algorithms, but these algorithms are usually not available as software tools for a broad range of users. Most implementations are platform specific and require a particular input format that is not used by other software tools. Changing a model to the required format is often a labor-intensive process. GINsim, a package designed for the analysis of gene regulatory networks, uses a method based on binary decision diagrams to efficiently analyze asynchronous logical models, and provides this method without requiring the user to understand the underlying algorithm, but GINsim is specific to logical models and cannot import other discrete modeling types or formats} \cite{GINsim}.
 
 \comment{In summary, only GINsim can non-heuristically analyze models with more than 32 variables for steady states and limit cycles when no initial state is given; and GINsim can analyze large models only when updates are assumed to be asynchronous}

Here, we present the web-based tool ADAM, Analysis of Dynamic
Algebraic Models \cite{ADAM}, a tool to study the dynamics of a wide
range of discrete models. \hl{ADAM provides efficient analysis methods based on mathematical algorithms as a web-based tool for several different input formats, and it makes analysis of complex models accessible to a larger community, as it is platform independent as a web-service and does not require understanding of the underlying mathematics.} ADAM is the successor to DVD, Discrete Visualizer of
Dynamics \cite{DVD}, a tool to visualize the temporal evolution of small
polynomial dynamical systems. 

\hl{Different types of discrete models, including  (probabilistic) Boolean networks, logical
networks, Petri nets, cellular automata, and agent-based 
models,} can be converted into the \hl{unifying framework of algebraic models, namely polynomial
dynamical systems} \cite{Alan:Bioinf2010, Hinkelmann:2010}. This allows to
apply tools from computational commutative algebra to analyze the models
more efficiently than by simulation and without using heuristic methods.
We used ADAM on several logical models and on published Boolean models with up to 60 variables \cite{GINsimRepo,AO,Pedicini-Tcell}. 
\hl{These models are too large for a straight-forward analysis by exhaustive
enumeration of the state space, and the corresponding publications contain lengthy explanations and supplementary material that outline the calculations and algorithms used to identify the attractors.} ADAM greatly simplifies the analysis; it
identified the steady states of these models in less than one second. \hl{We believe that complex discrete models will gain more popularity, if sophisticated analysis methods are easily accessible to modelers. }

In
addition to \hl{giving access to mathematical theory for efficient analysis, algebraic models provide} a unifying framework and systematic approach for several model types, which allows for an
effective comparison of heterogeneous models, such as a Boolean network model
and an agent-based model. For community integration to the biological
sciences, ADAM contains a model repository of previously published models
available in ADAM specific format \cite{ADAMRepo}. This allows new users to familiarize themselves quickly with ADAM and to validate and experiment with existing
models. In the following \hl{section}, we discuss general features of ADAM briefly and
explain new features in more detail.

\section*{Results and Discussion}
\hl{A wide range of software and algorithms exist to analyze discrete models. These tools are either limited by complexity as they rely on simulation as analysis method, or they are inaccessible to biologists not familiar with programming, as they are often times only available as platform dependent command-line tools. Furthermore, implementations require different input formats, which hampers the use of different software tools on the same model. 

ADAM is an analysis tool for discrete models, available as a web-based tool, which hides and encapsulates the mathematical algorithms from the user. Therefore, users who lack understanding of the underlying mathematics or programming expertise can use efficient algorithms to analyze complex models.
 
We propose a novel method to identify attractors of a discrete model. This method relies on the fact that many discrete models can be translated into the algebraic framework of polynomial dynamical systems. Using these polynomials, one can construct a system of polynomial equations, such that its solutions correspond to fixed points or limit cycles. Thus, the problem of identifying attractors becomes equivalent to solving a system of polynomial equations over a finite field. This is a long-studied problem in computer algebra, and can usually be solved efficiently by computing a Gr\"obner basis.

This method is not a new mathematical algorithm to solve polynomial dynamical systems, but a novel approach that uses the fact that attractors are the solutions of polynomial systems derived from the model when expressed in the algebraic framework. ADAM allows users unfamiliar with polynomial dynamical systems or Gr\"obner basis to benefit from this efficient algorithm. }

\subsection*{General Features}
ADAM automatically converts discrete models into polynomial dynamical systems,
that is time and state discrete dynamical systems described by polynomials
over a finite field (see Appendix \ref{sec:pds} for definition and example). The dynamics of the models are then analyzed by using various computational algebra
techniques. Even for large systems, ADAM computes key dynamic features, such
as steady states, in a matter of seconds. ADAM is available online and free of
charge. It is platform independent and does not require the installation of
software or a computer algebra system.

ADAM can analyze discrete models. \hl{It translates the following }{\bf inputs} \hl{into (probabilistic) polynomial dynamical systems and can then analyze all of them except models originating from Petri nets. }
\begin{itemize}
  \item Logical models generated with GINsim \cite{GINsim}
  \item Petri nets generated with Snoopy \cite{Snoopy}
  \item polynomial dynamical systems
  \item Boolean networks
  \item probabilistic polynomial dynamical systems, probabilistic Boolean networks (PBN) \cite{shmulevich}.
\end{itemize}
We plan to implement analysis methods for  Petri nets in future versions.

ADAM’s main application is the analysis of the dynamic features of a model, which includes the identification of stable attractors. These are either steady states, i.e.,
time-invariant states, or limit cycles, i.e., time-invariant sets of states. ADAM is capable of identifying all steady states and
limit cycles of length up to a user-specified length $m$. The process of
finding long limit cycles is quite slow for large models, however, in
biological models limit cycles are likely to be short, so that $m$ can be
chosen to be small in general, i.e., less than ten.

The temporal evolution of the model can be visualized by the {\it phase
space}, a graph of all possible states and their transitions, also called {\it state space} or {\it state transition graph}. For small enough
models, i.e., less than eleven variables, ADAM generates a graph of the complete phase space; for larger models, ADAM uses algebraic algorithms to determine dynamic properties. Independent of
network size, ADAM generates a {\it wiring diagram}. The wiring diagram, also
known as {\it dependency graph}, shows the static relationship between the
variables. All edges in ADAM’s wiring diagrams are functional edges, that is
there exists at least one state such that a change in the input variable causes a
change in the output variable (see Appendix \ref{sec:func} for more details).
This means that ADAM determines all non-functional edges, which is oftentimes
of interest.

With ADAM, one can also study the temporal evolution of user-specified initial states. The trajectory of a state describes the state's evolution, and it can be computed by repeatedly applying the transition function until an attractor is reached.  

All of these features can be computed assuming synchronous updates or sequential updates according to an update-schedule specified by the user. Note that the steady
states are the same independent of the update schedule. This is due to the
fact that updating any variable at a steady state does not change its value.
It is irrelevant for a steady state analysis whether updates are considered to happen sequentially or simultaneously.

For probabilistic networks, i.e., models in which each variable has several choices of local update rules, ADAM can generate a graph of all possible updates. This means that states in the phase space can have out-degree greater than one, since different transitions are possible. ADAM can find all true steady states, in the context of probabilistic networks, meaning all states that are time-invariant independent of the choice of update function. For further information of probabilistic networks, see \cite{shmulevich}.

For Boolean networks, ADAM calculates all functional circuits (see Appendix \ref{sec:func}). Positive functional circuits are a necessary condition for multi-stationarity. For a certain class of Boolean networks, namely conjunctive/disjunctive networks, ADAM computes a complete description of the phase space as described in \cite{conjunctive}. For further details on conjunctive networks, see Appendix \ref{sec:conj}.

In summary, ADAM can generate the following {\bf outputs}.
\begin{itemize}
\item wiring diagram
\item phase space for small models
\item steady states (for deterministic and probabilistic systems)
\item limit cycles of specified length $m$
\item trajectories originating from a given initial state until a stable
attractor is found
\item dynamics for synchronous or sequential updates
\item functional circuits for Boolean networks
\item a complete description of the phase space for conjunctive/disjunctive
networks.
\end{itemize}

\subsection*{Comparison to Other Systems}
In this section, we compare ADAM to other state-of-the-art systems.
Most software tools discussed here provide different functionality than ADAM does, thus 
a run-time comparison is not feasible. An overview of the
tools and their functionality is given in table \ref{table:comparison}, and they are discussed in detail below.

\begin{table}
	\begin{tabular}{|l|c|c|l|l|}
		\hline
	
		& \begin{sideways} Steady States \end{sideways} 
			& \begin{sideways} Limit Cycles \end{sideways} 
				& \multicolumn{1}{|c|}{\begin{sideways} Input\end{sideways}}
			& \multicolumn{1}{|c|}{\begin{sideways} Requirements\end{sideways}}  \\
		\hline
			\hline
	
		\multirow{3}{*}{ADAM} & \multirow{3}{*}{Yes$^\ddagger$}& \multirow{3}{*}{Yes$^\diamond$} & Boolean (or polynomial) functions &\multirow{3}{*}{None, web based} \\
			&&&Logical Models (GINsim) &\\
			&&&Petri Net (Snoopy) &\\
		\hline

		\multirow{2}{*}{GINsim}& \multirow{2}{*}{Yes$^\ddagger$} &
    \multirow{2}{*}{for asynchronous networks$^*$ }& parameters (non-zero truth tables) & \multirow{2}{*}{Java virtual machine$^\circ$}\\
			&&&Logical Model &\\
		\hline

		BoolNet R package & $\leq 29^ \dagger$  & $\leq 29^ \dagger$ & Boolean functions & R statistics software\\
		\hline

    \multirow{2}{*}{Snoopy} & \multicolumn{2}{|c|}{\multirow{2}{*}{reachable from a given initial state }} & SBML &\multirow{2}{*}{$^\circ$}  \\
      &\multicolumn{2}{|c|}{} & P/T graph &\\
    \hline
	
		DDLab & $\leq 31$ & $\leq 31$  & truth tables & $^\circ$\\
		\hline
		
		\hl{GenYsis-P } & \multicolumn{2}{|c|}{\multirow{2}{*}{Yes$^\otimes$}} & \multirow{2}{*}{GenYsis-P specific text file} & Linux that supports the \\ 
		Toolbox & \multicolumn{2}{|c|}{}  & &distributed binaries\\
		\hline
	
		BN/PBN Matlab  & \multirow{2}{*}{$\leq 27$} & \multirow{2}{*}{$\leq 27$ } & \multirow{2}{*}{truth tables} & \multirow{2}{*}{Matlab}\\
		Toolbox & & & \\
		\hline

	\end{tabular} 
	\caption{
	Comparison of different software tools regarding attractor analysis:
	$\ddagger$ less than 1 second
  		on published gene regulatory networks with up to 72 variables; 
	$\diamond$
  		only for short limit cycles;
	$ \dagger$ heuristic methods are available for larger networks; 
	$*$ asynchronous updates in the sense of logical models;
	$\otimes$ \hl{we had no platform available that supported the distributed binaries and could not conduct benchmark calculations};
		$\circ$ installation necessary, available for common operating systems.}
  \label{table:comparison}
\end{table}
	
{\it GINsim} (Gene Interaction Network simulation) is a package designed for the
analysis of gene regulatory networks \cite{GINsim}. As input, it accepts
logical models. Logical models are an extension to Boolean models; they
consist of similar switch-like rules, but allow for a finer discretization with more
than two states per variable, e.g., low, medium, and high. Logical models can be updated synchronously or asynchronously. For the latter, the temporal
evolution of a logical model is non-deterministic because the variables are updated randomly in an asynchronous
fashion. In either case, updates of every variable are continuous, meaning that no variable changes its value by more than one unit in one time-step, see section {\it Remarks about Logical Models} for a detailed discussion. 

GINsim
provides algorithms that use binary decision diagrams (BDD) for the
determination of steady states and oscillatory behavior \cite{Chaouiya}. 
For
synchronous updates, analysis of limit cycles is only possible by simulating every trajectory, i.e., generating the complete state space, called state transition graph in GINsim,
and therefore limited by network size.
We tested GINsim on logical models with up to $72$ variables; determining the
steady states took less than one second. More complex logical networks were not available to us.

Networks are entered
manually into GINsim, they cannot be imported from any other format.
Furthermore, models are specified by entering their parameters, i.e., entering
all values that result in a non-zero target value. Especially for large models, this can be a time consuming process.

{\it BoolNet R package} provides methods for inference and analysis of synchronous,
asynchronous, and probabilistic Boolean networks \cite{boolnet}. It is a
package for the free statistics software R, and it is run via the R command-line. It is helpful, if the user is already familiar with R. 
Steady state
analysis is implemented as exhaustive search of the state space, heuristic
search, random walk, or Markov Chain analysis \cite{shmulevich}. Non-heuristic analysis is limited to networks with 29 variables.
For
larger networks, steady states can be inferred heuristically, which does not guarantee that all steady states are identified.

{\it Snoopy} is a unifying Petri net framework, containing a family of Petri
net modeling tools and algorithms \cite{Snoopy}. Snoopy provides built-in simulation
and animation. Analysis of Petri nets can be performed, e.g., with the
tool {\it Charlie} \cite{Charlie}. Charlie identifies structural properties and
has algorithms for invariant based or reachability graph based analysis. \hl{Reachability graph-based} analysis for Petri nets usually depends on a given initial state and does not provide a complete picture of possible dynamics for other initial markings. In addition to marking-dependent analysis, Charlie uses algorithms based on linear algebra to predict the dynamic properties independent of markings, such as T and P-invariants. ADAM converts a Petri net to a collection of polynomial dynamical systems, one system for each transition. Analysis of such a non-deterministic system is currently not implemented in ADAM, and we therefore do not list any run-time results for Petri nets. 

{\it DDLab} is an interactive graphics software for discrete models, including
cellular automata, Boolean and multi-valued networks \cite{DDLab}. As it is
mainly a visualization tool, analysis is based on exhaustive enumeration of
the state space, and model size is limited to 31 variables.  

{\it GenYsis-P Toolbox} \hl{is a command-line tool currently available only for Linux to analyze (probabilistic) gene regulatory networks. Algorithms use (reduced order) binary decision diagrams. As analysis methods are not based on exhaustive enumeration, GenYsis-P can analyze large networks. Unfortunately, we did not have access to a platform that supports the distributed binaries, and source code was not available. }

{\it BN/PBN Toolbox} is a toolbox written in Matlab \cite{shmulevich-matlab}.
It uses the state transition matrix to compute attractors. Statistics for
networks with more than 27 variables cannot be computed (``Maximum variable
size allowed by the program is exceeded"). In addition to analyzing
deterministic Boolean networks, the toolbox can analyze probabilistic Boolean
networks and calculate statistics such as numbers and sizes of attractors,
basins, transient lengths, Derrida curves, percolation on 2-D lattices, and
influence matrices. 

ADAM calculates all steady states of networks
non-heuristically, by applying algebraic algorithms, see section {\it Methods}. All of the above software
tools provide other functionality that ADAM currently lacks, \hl{but for the analysis of synchronous networks, they all are either restricted to less than 32 variables or require familiarity with programming.} GINsim calculates steady states for large networks, but models must
be specified by their parameters and cannot be entered in a compact form as such
(Boolean) functions.

Several Boolean models too large for analysis by exhaustive enumeration have been published, 
for example on the expression of the segment polarity genes in Drosophila
melanogaster or T-cell regulation
\cite{AO,Pedicini-Tcell}. ADAM identified all steady states and small limit cycles
for these systems
in less than one second, whereas the publications that comprise the Boolean models, contain lengthy explanations outlining the analysis. The model files in ADAM format can be accessed at
\cite{ADAMRepo}.

\subsubsection*{Remarks about Logical Models} 
ADAM allows for synchronous or sequential updates \hl{according to a given update schedule}. In models with synchronous updates, all variables are updated simultaneously at every time step. In models with sequential updates, all variables are updated at every time-step, but in the order of the given update schedule. \hl{Models with sequential updates can be converted into synchronous models with identical state space}. In models with asynchronous updates, \hl{as it is common} for logical models, one variable is updated at random at every time step, \hl{which results in a non-deterministic model.}
Sequential and asynchronous updates of the same system result in different dynamics.

In GINsim, all models are \hl{``continuous'' in the sense that at each time-step, each variable increases or decreases by at most one unit. Though logical models are discrete, there are no jumps skipping intermediate states. For example, in a model with three states, low, medium, and high, no variable can drop from high to low in a single update step. This interpretation is different from the common meaning of continuous, which usually refers to models of ordinary or partial differential equations}. The parameters entered in GINsim specify the target value towards which the variable changes, i.e., the value increases by one, decreases by one, or remains constant if the target value is larger, smaller, or equal than the initial value, respectively. The phase space generated with ADAM might differ from the state transition graph generated in GINsim. To obtain the exact same phase space, every variable in the logical model must contain an explicit self-loop, and all parameters must be entered such that the target value differs by at most one from the value of the variable to be updated. Any logical model can be specified in this way without changing its state transition graph. Boolean models are always continuous.   

In multi-valued logical models, variables can have different maximum values. In an algebraic model, all variables are defined over the same algebraic field, i.e., have the same maximum value. When a multi-valued logical model is translated into an algebraic model, extraneous states might be introduced such that all variables are defined over the same field. An example of such an extension is given in table \ref{table:extraStates}, the extra states are the states in the last row, which are given the same values as the states above to extend the model in a meaningful way. The extra states should be ignored when analyzing the dynamics. For more details, see \cite{Alan:Bioinf2010}.  

\comment{\begin{table}
	\begin{tabular}{||c|c||c|c||c|c||c|c||}
		\hline
		\multicolumn{2}{||c||}{$x_2$ }& $x_2=$low & $x_2=0$ & $x_2=$medium & $x_2=1$ &$x_2=$high & $x_2=2$ \\

		\hline \hline
			$x_1$ absent & $x_1 =0$ & $x_2=$low & $x_2=0$ & $x_2=$medium & $x_2=1$ &$x_2=$high & $x_2=2$ \\
		\hline 
			$x_1$ present & $x_1 =1$ & $x_2=$medium & $x_2=1$ & $x_2=$high & $x_2=2$ &$x_2=$high & $x_2=2$ \\
		
		\hline
			& $x_1 =2$ & & $x_2=1$ & & $x_2=2$ & & $x_2=2$ \\
		\hline 
	\end{tabular}
	\caption{Updates for variable $x_2$ in a logical model. The states $0$ and $1$ represent absent and present for the Boolean variable $x_1$; $0$, $1$, and $2$ represent low, medium, and high for the multi-valued variable $x_2$. The last row is introduced in the polynomial dynamical system such that all variables are defined over $\mathbb F_3$. The extra states $(2,0), (2,1), (2,2)$ in the state space should be ignored when interpreting the dynamics.}
	\label{table:extraStates}
\end{table}
}
\begin{table}
	\begin{tabular}{|c||c|c|c|}
		\hline
		 		next state of $x_2$& low $x_2$ & medium $x_2$  & high $x_2$  \\
	
		\hline \hline
			$x_1$ absent & low $x_2$ & medium $x_2$ & high $x_2$ \\
		\hline 
			$x_1$ present & medium $x_2$ & high $x_2$  & high $x_2$ \\
		\hline
			extension $x_1$ present & medium $x_2$ & high $x_2$  & high $x_2$ \\
		\hline 
	\end{tabular}
	\caption{\hl{Updates for variable $x_2$ in a logical model, where $x_2$ depends on $x_1$ and itself. The states $0$ and $1$ represent absent and present for the Boolean variable $x_1$; $0$, $1$, and $2$ represent low, medium, and high for the multi-valued variable $x_2$. The last row is introduced in the polynomial dynamical system such that all variables are defined over $\mathbb F_3$. The extra states $(2,0), (2,1), (2,2)$ in the state space should be ignored when interpreting the dynamics.}}	
	\label{table:extraStates}
\end{table}

\subsubsection*{Remarks about Petri Nets}
In the Petri net community, {\it state space} usually refers to all states (markings) reachable for a given initial state. In this manuscript, {\it state space} refers to the set of all possible states, independent of an initial state.

Translating a bounded Petri net to an algebraic model results in a set of polynomial dynamical systems, where every transition corresponds to one system. In a Petri net, different firing sequences can lead to different markings; a firing sequence relates to the order in which the different systems are iterated. The {\it update schedule} is not related to the firing sequence. As firing does not consume any time, polynomial systems describing a Petri net always use synchronous updates \cite{Heiner:2008}. 

The term {\it functional edge} is not related to the concept of {\it liveness of a transition}. The liveness of a transition depends on an initial marking. An edge, connecting two places (source and target), is functional, if there exists a marking, such that changing only the marking of the source place, changes the marking of the target place (see Appendix \ref{sec:func} for more details on functional edges). 

\subsection*{Application} \label{benchmarks}
We show how to use ADAM on a well-understood model of the expression
pattern of the segment polarity genes in Drosophila melanogaster. Albert and Othmer developed a model for embryonic pattern formation in the fruit fly Drosophila melanogaster \cite{AO}. Their Boolean
model consists of 60 variables, resulting in a state
space with more than $10^{18}$ states. They analyze the model for steady states by manually solving a system of Boolean equations. They also analyze the temporal evolution of a specific initial state corresponding to the wild type expression pattern by repeatedly applying the Boolean update rules until a steady state is found. The update schedule of the model is synchronous  with the exception of activation of SMO and the binding of PTC to HH (activation of PH), which are assumed to happen instantaneously. This can be accounted for by substituting the equations for SMO and PH into the update rules for other genes and proteins, rather than using SMO and PH themselves.

To analyze the model, we first rename the variables in the Boolean rules given in \cite{AO} such as $\text{wg}_i$ or $\text{SLP}_i$ to \hl{$x_1 \ldots x_{60}$, to standardize their format}. Then we use ADAM: the model type is {\it Polynomial Dynamical Systems}, the number of states in a Boolean model is $2$, representing present or absent. One can choose {\it Boolean}, and enter the Boolean rules in the text-area or upload a text file with the Boolean rules. Alternatively, one can first convert the Boolean rules to polynomials over $\mathbb F_2$, and enter the polynomials with the choice {\it Polynomial}.
The file with the polynomial equations for the model can be accessed at \cite{ADAMRepo}. 

The rules in the model file are specified in {\it Polynomial} form. Once the polynomials are uploaded, we need to set the {\it Analysis} type. The model with $60$ variables is too complex for exhaustive enumeration, and we choose {\it Algorithm}. This means that instead of exhaustive enumeration of the state space, analysis of the dynamics is done via computer algebra by solving systems of equations.
In {\it Options}, we set {\it Limit cycle length} to one because we are interested in the steady states, i.e., time-invariant states.
We chose {\it synchronous} as updating scheme. Once these choices have been made,
we obtain the steady states by clicking {\it Analyze}.
ADAM returns a link to the {\it wiring diagram} or dependency graph, which captures the static relations between the different variables. In addition, ADAM returns the number of steady states and the steady states themselves, see figure \ref{fig:alg}. \hl{These steady states are identical to those found in} \cite{AO}, \hl{half of which have been observed experimentally.} 
\begin{figure}[htb]
	\centering
	\includegraphics[width=0.95\textwidth]{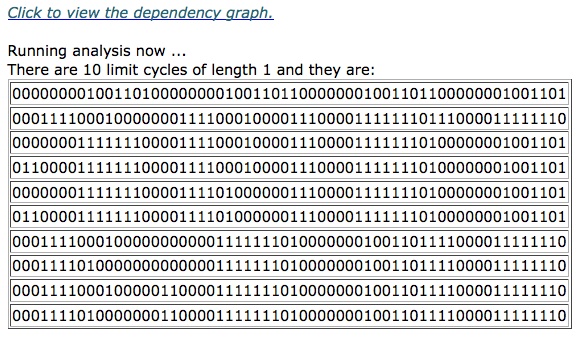}
	\caption{\hl{ADAM: Analysis of steady states of Drosophila model. Each row in the table corresponds to a stable attractor. Attractors are written as binary strings, where $0$ represents non-expression of a gene (or low concentration of a protein), and $1$ expression (or high concentration)}}
	\label{fig:alg}
\end{figure}

Each row in the table corresponds to a stable attractor. Attractors are written as binary strings, \hl{where $0$ represents non-expression of a gene (or low concentration of a protein), and $1$ expression (or high concentration), e.g.,} 
\begin{align}\label{equ:ss}
	0 0 0 1 1 1 1 0 0 0 1 0 0 0 0
	0 0 0 0 0 0 0 1 1 1 1 1 1 1 0
	1 0 0 0 0 0 0 0 1 0 0 1 1 0 1
	1 1 1 0 0 0 0 1 1 1 1 1 1 1 0
\end{align}

\hl{corresponds to the genes (and proteins) being expressed ( or present in high concentration) in four cells form anterior to posterior compartments (compartment $1$ to $4$) as shown in table} \ref{table:ss}.
\begin{table}
	\begin{tabular}{|l|l|}
		\hline
		compartment 1 & en,	EN,	hh,	HH,	SMO\\
		\hline
		compartment 3 &ptc, PTC, PH, SMO, ci, CI, CIA\\
		\hline
		compartment 3 &SLP, PTC, ci, CI, CIR\\
		\hline
		compartment 4 &SLP, wg, WG, ptc, PTC, PH, SMO, ci, CI, CIA\\
		\hline
	\end{tabular}\caption{\hl{Genes and proteins present in steady state} \ref{equ:ss}}
	\label{table:ss}
\end{table}
\hl{This is the steady state obtained in} \cite{AO} \hl{when starting the system with an initial state representing the experimental observations of stage 8 embryos.} \comment{Note that we include variables for SLP, which are not shown in 
as they are fixed to be OFF in the anterior part, and ON in the posterior part throughout the paper, i.e., $0011$.}

ADAM can also generate trajectories for a given initial state. For example, we
can choose the initial state that was used in \cite{AO} representing stage 8 embryos. Again, we
enter {\it Polynomial Dynamical Systems} with $2$ as the number of states and upload the polynomials
describing the model. Instead of {\it Algorithms}, we now choose {\it
Simulation}. Since we are not interested in the number of steady states or the
complete phase space, but in a single trajectory originating from a specific
initial state, we choose {\it One trajectory starting at an initial state} as the simulation option. As
initial state we enter
\begin{align*}
0 0 0 1 0 1 0 0 0 0 0 0 0 0 0
0 0 0 0 0 0 0 1 0 0 0 1 0 0 0
1 0 0 0 0 0 0 1 0 0 0 1 0 0 0
1 1 0 0 0 0 0 1 0 0 0 1 0 0 0.
\end{align*}
By clicking {\it Analyze}, we obtain the temporal evolution of this particular
state until it reaches a steady state. As predicted in \cite{AO}, the
steady state is the state described in table \ref{table:ss}, see Fig.
\ref{fig:traj}.
\begin{figure}[htb]
\centering
\includegraphics[width=0.95\textwidth]{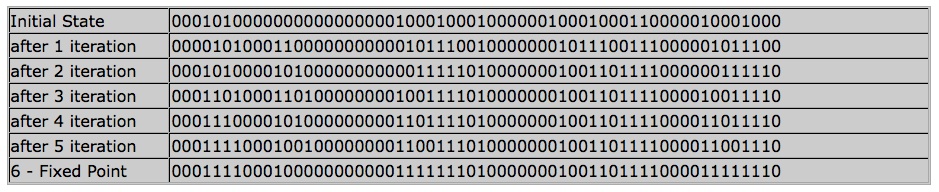}
\caption{ADAM: Trajectory of Drosophila model}
\label{fig:traj}
\end{figure}

To summarize, ADAM correctly identified the steady states
in less than one second. All steady states have been determined previously in \cite{AO} by labor-intensive manual investigation of the system. \comment{In \cite{AO}, the model is formulated as a set of Boolean rules. In order to determine the steady states, the system of Boolean expressions was solved manually.}

Furthermore, we used ADAM to verify that there are no limit
cycles of length two or three. The model has not been analyzed previously for
limit cycles. The absence of two- and three cycles strengthens confidence in
the model, since oscillatory behavior has not been observed experimentally. \hl{Computations for limit cycles of length greater than three have not been conducted, as composing the system several times with itself is computationally complex.}
The model file in ADAM format can be accessed at \cite{ADAMRepo}.


\subsubsection*{Benchmark Calculations}
We analyzed logical models
available in the GINsim model repository \cite{GINsimRepo} as of August 2010. The
repository consists of models in GINsim XML format previously published in
peer-reviewed journals. We converted all but two models into polynomial
dynamics systems. For these 26 models we computed the steady states. All
calculations finished in less than 1,5 seconds, see
Figure \ref{fig:chart}.

In addition to the published models in \cite{GINsimRepo}, we analyzed
randomly generated networks
that have the same structure that we
expect from biological systems, namely sparse, i.e., while the number of nodes in a biological network may be quite large, each node is affected only by a small number of other nodes, and robust, i.e., small number of attractors.. We tested a total of 50 networks with
50-150 nodes ($10^{15} - 10^{45}$ states) and \hl{an average of average in-degrees of $1.6848$}. The
steady state calculations took less than half a second for each network on
a 2.7 GHz computer.
\begin{figure}[htb]
\centering
\includegraphics[width=0.7\textwidth]{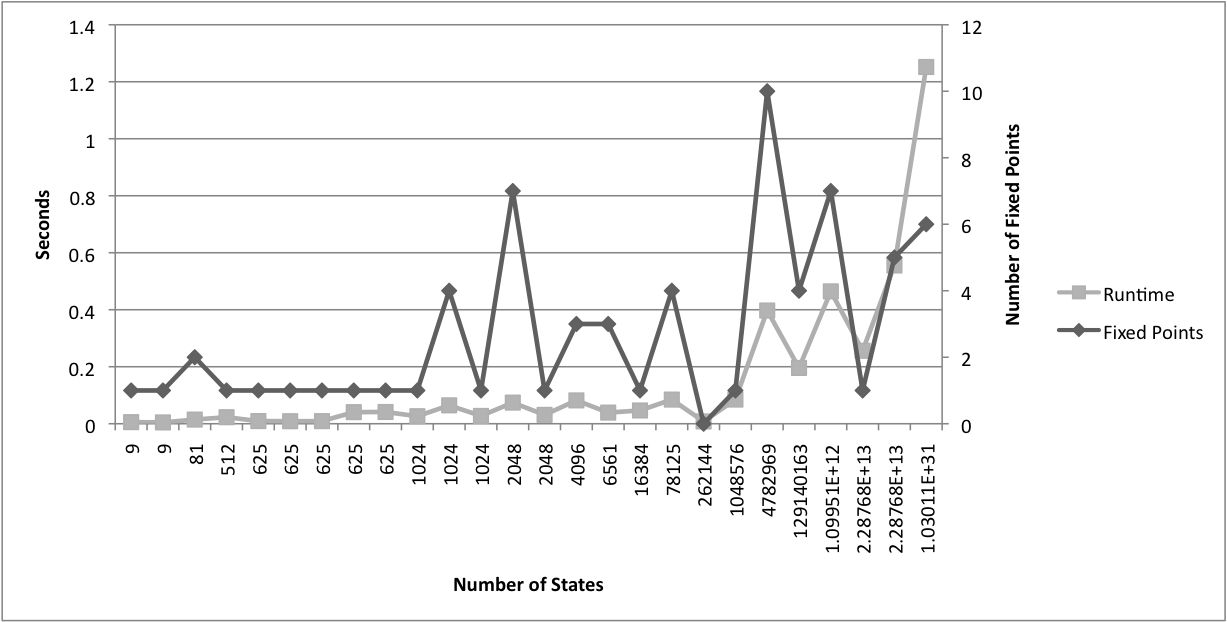}
\caption{\hl{Runtime of steady state calculations of several logical models from}
\cite{GINsimRepo}. \hl{Executed on a 2.7 GHz computer.}}
\label{fig:chart}
\end{figure}

\subsection*{Architecture}
ADAM is available as an web-based tool, that does not require any software installation. ADAM's user interface is implemented in HTML. We use JavaScript to generate a dynamic website that adapts as the user makes various choices. This simplifies the process of entering a model. For example, after defining the model type, i.e., Polynomial Dynamical System, Probabilistic Network, Petri net, and Logical Model the next line changes to the number of states, $k$-bound, or nothing, appropriately. Input can be entered directly into the text-area, or uploaded as a text document.

All mathematical algorithms are programmed in Macaulay2 \cite{M2}. Macaulay2 is a powerful computer algebra system. The routines for which fast execution is crucial are implemented in C/C++ as part of the Macaulay2 core. Logical Models and Petri nets in XML format are parsed using Ruby's XmlSimple library. The interplay between HTML and Macaulay2 is also programmed in Ruby.

Output graphs are generated with Graphviz's {\it dot} command. When {\it Simulation} is chosen as analysis method, Graphviz's {\it ccomps - connected components filter for graphs} is used to count the connected components. A Perl script directs the execution of the Graphviz commands.
\subsection*{Model Repository}
A model repository is part of the ADAM website \cite{ADAMRepo}. The repository consists of a collection of several previously published models in ADAM format. The models are extracted from publications, and rewritten in ADAM specific format, i.e., all variables are renamed to $x_i$ and the update rules from the original publication are reformulated as Boolean rules or polynomials. The central repository with models in a unified framework allows for quick verification and experimentation with published models. By changing parameters or initial states, users can gain a better understanding of the models.

New users can also use the repository to quickly familiarize themselves with the main functionalities of ADAM. In addition to the model itself, the database entries contain a short summary of the biological system and relevant graphs, together with an analysis of dynamic features determined by ADAM and their biological explanation. The repository is work in progress by researchers from several institutions generating more entries for the repository. We invite all interested researchers to submit their models.

Because of their intuitive nature, discrete models are an excellent introduction to mathematical modeling for students of the life sciences. ADAM's model repository is a great starting point to familiarize students with the abstraction of discrete models such as Boolean networks.
 
%
%

\section*{Conclusions}
Discrete modeling techniques are a useful tool for analyzing complex biological systems \hl{and there is a need in the biological community for easily accessible analysis software. ADAM provides efficient methods as a web-based tool and will allow a larger community to use complex modeling techniques, as it is platform independent and does not require the user to understand the underlying mathematics.}
Upon translating discrete models, such as logical networks,
Petri nets, or agent-based models into algebraic models, rich mathematical
theory becomes available for model analysis.

After extensive experimentation with both discrete models arising in systems biology and randomly generated networks, we found that the algebraic algorithms presented in this manuscript
are fast for sparse systems \hl{with few attractors}, a structure maintained by most biological
systems. All algorithms have been included in the software package ADAM\cite{ADAM},
which is user-friendly and available as a free web-based tool.
ADAM is highly suitable to be used in a classroom as a first
introduction to discrete models because students can use it without going through an installation process.

\comment{Several software tools for discrete models exists, all of them specializing in a single discrete model type.
{\it GINsim}, a tool for modeling and simulation of logical models can quickly identify steady states \cite{GINsim}, but it has no other means than simulation of every trajectory to identify limit cycles for synchronous networks. {\it BoolNet R package}, a package for inference and analysis of synchronous, asynchronous, and probabilistic Boolean networks, does a steady state analysis by exhaustive enumeration of the state space or using heuristic methods \cite{boolnet}. Analysis is limited by model size (exhaustive enumeration or Markov Chain analysis ) or restricted to heuristic methods, which might fail to detect some key dynamic features.
{\it Snoopy} and {\it Charlie}, tools for Petri nets, provide marking-independent analysis based on linear algebra, e.g., methods to determine T-invariants, which can be interpreted as the relative firing rates, occurring permanently, and this activity level corresponds to the steady state behavior.
{\it DDLab}, an interactive graphics software for cellular automata, Boolean and multi-valued networks, does not provide analysis methods other than through visualization, which is limited by model size or restricted to a small portion of the state space.}

ADAM provides methods to analyze the key dynamic features, such as steady states and limit cycles, for large-scale (probabilistic) Boolean networks and logical models. ADAM unifies different modeling types by providing analysis methods for all of them and thus can be used by a larger community.

We hope to expand ADAM to a more comprehensive Discrete Toolkit which incorporates new analysis methods, better visualization, and automatic conversion for more model types.
We also hope to analyze controlled algebraic models and expand theory to stochastic systems.

\section*{Methods}
Logical models, Petri nets, and Boolean networks are converted automatically
into the corresponding polynomial dynamical system as described in
\cite{Alan:Bioinf2010}, so that algorithms from computational
algebra can be used to analyze the dynamics. \hl{In polynomial dynamical systems over a finite field, states of a variable are assigned to values in the field, and the update (or transition) rule for each variable is given as a polynomial rather than a Boolean or logical expression. For more details, see appendix} \ref{sec:pds}. \hl{Using these polynomials, one can construct systems of polynomial equations, such that their solutions correspond to fixed points or limit cycles. Thus, the problem of identifying attractors becomes equivalent to solving a system of polynomial equations over a finite field. This is a long-studied problem in computer algebra, and can usually be solved efficiently by computing a Gr\"obner basis.}

Gr\"obner basis calculation is for polynomial systems what
Gauss-Jordan elimination is for linear systems: a structured way to transform
the original system to triangular shape without changing its solution space.
The triangular shape of the resulting systems allows for stepwise retrieval of the solutions of the system. For a more in depth discussion of Gr\"obner bases, see for example \cite{IVA}.

In the worst case, computing Gr\"obner bases for a set of polynomials has
complexity doubly exponential in the number of solutions to the system.
However, in practice, Gr\"{o}bner bases are computable in a reasonable time. It has been suggested, that in robust gene
regulatory networks genes are regulated by only a handful of regulators
\cite{Leclerc:2008}. Thus, the polynomial dynamical systems representing such biological networks are
sparse, i.e., each function depends only on a small subset of the model variables. From our experience, a Gr\"obner basis calculation 
for sparse systems with few attractors, a structure common for
biological systems, is actually quite fast. 

\appendix
\section{Mathematical Background}
\subsection{Polynomial Dynamical Systems}\label{sec:pds}
To be self-contained, we briefly explain polynomial dynamical systems and their key features.

A {\bf polynomial dynamical system} (PDS) \cite{JLSS} over a finite field $k$ is a function
$$f = (f_1, \ldots, f_n) : k^n \rightarrow k^n,$$
with coordinate functions $f_i \in k[x_1, \ldots , x_n]$. Iteration of $f$ results
in a time-discrete dynamical system. \hl{A PDS can be used to describe the dynamic behavior of a biological system: every variable $x_i$ corresponds to a biological substrate, for example a protein or gene, and the polynomials $f_i$ describe the evolution of $x_i$ depending on the previous state of the variables $x_1, \ldots, x_n$.}

\comment{PDS are special cases of finite
dynamical systems, which are maps $X^n \rightarrow   X^n$ over arbitrary
finite sets $X$.}
PDS have several dynamic features of biological
relevance. These include the number of components, component sizes, steady states, limit cycles, and limit cycle lengths.
\begin{example}
Let $k= \mathbb F_2$ and $f = (f_1, f_2, f_3) : \mathbb F_2^3 \rightarrow
\mathbb F_2^3$ with
\begin{align*}
f_1 &= x_1x_2x_3+x_1x_2+x_2x_3+x_2 \\
f_2 &= x_1x_2x_3+x_1x_2+x_1x_3+x_1+x_2 \\
f_3 &= x_1x_2x_3+x_1x_3+x_2x_3+x_1+x_2.
\end{align*}
The wiring diagram of $f$, which shows the static interaction of the three
variables, is
depicted in Figure \ref{fig:ex} (left) along with its phase space in Figure
~\ref{fig:ex} (right).
The phase space shows the temporal evolution of the system. Each state is
represented as a vector of the values of the three variables $(x_1, x_2,
x_3)$.
The PDS described by $f$ has
two stable attractors: a steady state, $(000)$, and a limit cycle of length
three, consisting of the states $(010)$, $(111)$, and $(011)$.
\end{example}
\begin{figure}[ht]
\centering
\includegraphics[scale=0.55]{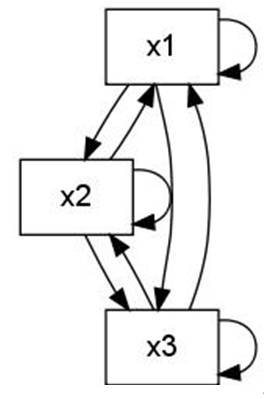}
\includegraphics[scale=0.55]{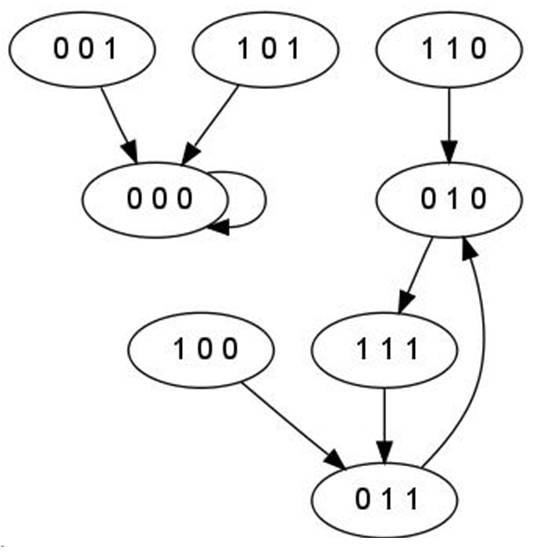}
\caption{(left)
Wiring diagram: static relationship between variables
(right)
Phase space: temporal evolution of the system
}
\label{fig:ex}
\end{figure}

A {\bf probabilistic PDS} over a finite field $k$ is a collection of functions
$$f = (\{f_{1,1}, \ldots, f_{1, r_1}\}, \ldots, \{f_{n, 1}, \ldots, f_{n, r_n}
\}) : k^n \rightarrow k^n,$$
together with a probability distribution for every coordinate that assigns the
probability that a specific function is chosen to update that coordinate.
The coordinate functions $f_{i,j}$ are elements in $k[x_1, \ldots , x_n]$.
Probabilistic PDS, specifically Boolean probabilistic networks (PBN), have been studied
extensively in \cite{shmulevich}.
ADAM analyzes probabilistic PDS. It can simulate the
complete phase space for small enough models, by generating every possible
transition and labeling the edge with its probability according to the
distribution. If no distribution is given, ADAM assumes a uniform distribution
on all functions. For large networks, ADAM's {\it Algorithm} choice computes
steady states of probabilistic networks.
\subsection{Functional Edges} \label{sec:func}
An edge in the wiring diagram from $x_i$ to $x_j$ is considered
functional, if there exists a state $\hat x = (\hat x_1,  \ldots, \hat x_n)$ such
that $f_j( \hat x_1,  \ldots, a, \ldots \hat x_n) \neq f_j(\hat x_1, \ldots, b, \ldots
\hat x_n)$, where $a$ and $b$ are values for $x_i$, in other words, if there
is at least one state, such that changing only $x_i$ but keeping all other
values fixed, changes the next state of $x_j$.
In ADAM, all edges in the wiring diagram are functional.
For Boolean networks, ADAM identifies all functional elementary circuits. An elementary circuit is a finite closed paths in the wiring diagram where all the nodes are distinct. Functional Circuits are a necessary condition for multi-stationarity and limit cycles. For a further discussion of
functional circuits, see \cite{Chaouiya}. For multivalued networks, circuit analysis has not yet been implemented.
\section{Algorithms}
\subsection{Analysis of stable attractors}
Every attractor in a PDS is either a
steady state or a limit cycle. For small models, ADAM determines the complete
phase space by enumeration, for large models, ADAM computes steady states and
limit cycles of a given length.
A state is a steady state, if it transitions to itself after one update of the
system. A state is part of a limit cycle of length $m$, if,
after $m$ updates, it results in itself. Any steady state of a PDS satisfies
the equation $f(x) = x$, as no coordinate of $x$ is changing as it is updated.
Similarly, states of a
limit cycle of length $m$ satisfy the equation $f^m(x) = x$. ADAM computes all
steady states by solving the system $f_i(x) - x_i = 0$ for $i \in \{1, \ldots,
n\}$ simultaneously. To efficiently solve the resulting systems of polynomial
equations, we first compute the Gr\"obner
basis in lexicographic order for the ideal generated by the equations.
By the elimination and extension theorem, choosing a lexicographic order
allows to easily obtain the solutions \cite{IVA}.
We use the Gr\"obner basis algorithms distributed with Macaulay2, a
computer algebra system, and found that for quotient rings over a finite field
the implementation `Sugarless' is more efficient than the default algorithm
with `Sugar' \cite{M2,Sugar:1991}.
For limit cycles of length $m$, the solutions of $f^m(x)=x$ are found and then
grouped into cycles, by applying $f$ to each of the solutions.

\begin{example}
	\hl{Fixed points of the system shown in the example in} \ref{sec:pds} \hl{are solutions in $\mathbb F_2^3$ of the system $f(x)=x$:}
	\begin{align*}
	 x_1x_2x_3+x_1x_2+x_2x_3+x_2  &= x_1\\
	 x_1x_2x_3+x_1x_2+x_1x_3+x_1+x_2 &= x_2 \\
	 x_1x_2x_3+x_1x_3+x_2x_3+x_1+x_2 &= x_3.
	\end{align*}
	\hl{The only solution to this systems is the point $(x_1, x_2,x_3 ) = (0,0,0)$. This is in accordance with the state space depicted in figure} \ref{fig:ex}: \hl{$(0,0,0)$ is the only fixed point. To investigate limit cycles of length two, one has to look at the system $f^2(x) = x$,}
	\begin{align*}
	g_1(x) = f_1(f_1(x), f_2(x), f_3(x) ) &= x_1*x_2+x_2*x_3 = x_1 \\
	g_2(x) = f_2(f_1(x), f_2(x), f_3(x) ) &= x_1*x_2*x_3+x_1*x_2+x_1*x_3+x_1+x_2 = x_2\\
	g_3(x) = f_3(f_1(x), f_2(x), f_3(x) ) &= x_1*x_2*x_3+x_2 = x_3.
	\end{align*}
	\hl{Again, $(0,0,0)$ is the only solution, which means that there are no limit cycles of length two. 

	Investigating $f^3(x) = x$,}
	\begin{align*}
	f_1(g_1(x), g_2(x), g_3(x) ) &= x_1 \\
	f_2(g_1(x), g_2(x), g_3(x) ) &= x_2\\
	f_3(g_1(x), g_2(x), g_3(x) ) &= x_3, 
	\end{align*}
	\hl{results in the solutions $(0, 0, 0), (0, 1, 0), (0, 1, 1), (1, 1, 1)$. $(0,0,0)$ is a fixed point, and $(0, 1, 0), (0, 1, 1), (1, 1, 1)$ are elements of a limit cycle of length three. For all $m > 3$, $f^m(x)=x$ has no solutions, that means the system $f$ has exactly two attractors, a fixed point a a limit cycle of length three.} 
\end{example}
\subsection{Conjunctive/Disjunctive Networks} \label{sec:conj}
Some classes of networks have a certain structure that can be
exploited to achieve faster calculations. Jarrah et al.
show that for conjunctive (disjunctive) networks key dynamic features can be found with
almost no computational effort \cite{conjunctive}. Conjunctive (respectively disjunctive) networks consist of
functions using only the AND (respectively OR) operator.
ADAM comes with an implementation of this algorithm to analyze
dynamics in the case of conjunctive (disjunctive) networks. Currently,
this option is only implemented for networks with strongly connected dependency graphs.


\section*{Authors' contributions}
FH led the algorithm and software development group. BG, MB, and RM implemented the user interface and attractor analysis, executed benchmarking calculations, and drafted initial manuscript under FH's leadership. AV implemented the translation for logical algorithms to PDS used by ADAM. GB participated in the software design effort and
algorithm development. RL  conceived of the study, provided overall leadership of the project, and secured
funding for it. He also contributed to the writing and editing of the manuscript. 
All authors read and approved the final manuscript.

\section*{Acknowledgements}
  \ifthenelse{\boolean{publ}}{\small}{}
 Dimitrova, Clemson University; J. Adeyeye, Winston-Salem State University; B. Stigler, Southern Methodist University; R. Isokpehi, Jackson State University are currently expanding ADAM’s Model Repository.
Funding for this work was provided through U.S. Army Research Office Grant Nr. W911NF-09-1-0538,
National Science Foundation Grant Nr. CMMI-0908201, and National Science Foundation Grant Nr. 0755322.


{\ifthenelse{\boolean{publ}}{\footnotesize}{\small}
 \bibliographystyle{bmc_article}  
  \bibliography{ADAM} }     


\ifthenelse{\boolean{publ}}{\end{multicols}}{}



\section*{Figures}
  \subsection*{Figure 1 - ADAM: Analysis of steady states of Drosophila model. Each row in the table corresponds to a stable attractor. Attractors are written as binary strings, where 0 represents non-expression of a gene (or low concentration of a protein), and 1 expression (or high concentration)}
      Steady states of Drosophila Melanogaster as found with ADAM.

  \subsection*{Figure 2 - ADAM: Trajectory of Drosophila model}
      Temporal evolution of given initial state until steady state is reached.

	\subsection*{Figure 3 - Runtime of steady state calculations of several logical models from [17]. Executed on a 2.7 GHz computer.}
	      Runtime of steady state calculations of several logical models from [17]. Executed on a 2.7 GHz computer.
	
		\subsection*{Figure 4 - (left) Wiring diagram: static relationship between variables (right) Phase space: temporal evolution of the system.}
		(left) Wiring diagram: static relationship between variables (right) Phase space: temporal evolution of the system.
		
\subsection*{Table 1 -
	Comparison of different software tools regarding attractor analysis:
	$\ddagger$ less than 1 second
  		on published gene regulatory networks with up to 72 variables; 
	$\diamond$
  		only for short limit cycles;
	$ \dagger$ heuristic methods are available for larger networks; 
	$*$ asynchronous updates in the sense of logical models;
	$\otimes$ \hl{we had no platform available that supported the distributed binaries and could not conduct benchmark calculations};
		$\circ$ installation necessary, available for common operating systems.}
		Comparison of different softwares regarding attractor analysis.

\subsection*{Table 2 - 
	Updates for variable x2 in a logical model, where x2 depends on x1 and itself. The states 0 and 1 represent absent and present for the Boolean variable x1; 0, 1, and 2 represent low, medium, and high for the multi-valued variable x2. The last row is introduced in the polynomial dynamical system such that all variables are defined over F3. The extra states (2, 0), (2, 1), (2, 2) in the state space should be ignored when interpreting the dynamics.}
	Extending state space when converting logical model to polynomial dynamical system.
	
\subsection*{Table 3 -
 	Genes and proteins present in steady state}
	Steady state

\end{bmcformat}
\end{document}